\documentclass[11pt]{article}

\usepackage[dvips]{graphicx}

\newtheorem{definition}{Definition}
\newtheorem{theorem}{Theorem}

\newtheorem{remark}{Remark}

\newcommand{\bprf}{\bf Proof: \nopagebreak  \rm}
\newcommand{\eprf}{\nopagebreak\hspace*{1em} \hfill
q.e.d.\noindent\vspace{2ex}\\\phantom{}} 

\title{On Coercivity and the Frequency Domain Condition in Indefinite LQ-Control}

\author{Tobias Damm\footnote{
{\tt damm@mathematik.uni-kl.de}, University of Kaiserslautern,
Department of Mathematics, Gottlieb-Daimler-Stra\ss{}e, and Fraunhofer
ITWM, 67663 Kaiserslautern, Germany} 
\and Birgit Jacob\footnote{{\tt bjacob@uni-wuppertal.de}, University of Wuppertal, IMACM,
School of Mathematics and Natural Sciences,
Gau\ss{}stra\ss{}e 20, 42119 Wuppertal, Germany}}
\date{~}
\begin{document}

\maketitle

\pagestyle{myheadings}
\markboth{T. Damm, B. Jacob
}{Coercivity and Frequency Domain Condition}

\bigskip

\begin{abstract}
We introduce a coercivity condition as a time domain analogue of the
frequency criterion provided by the famous Kalman-Yakubovich-Popov
lemma. For a simple stochastic linear quadratic control problem we show how the
coercivity condition characterizes the solvability of Riccati equations.


\end{abstract}


\section{Introduction}
\label{sec:introduction}
Since the formulation of the Kalman-Yakubovich-Popov-lemma in the 1960s
the interplay of time domain and frequency domain methods has always
been fruitful and appealing in linear control theory. For the
linear-quadratic control problem and the algebraic Riccati equation
this has been worked out to a large extent already in \cite{Will71}.

However, the applicability of frequency domain methods is mostly
limited to linear time-invariant deterministic models. In the
consideration of time-varying or stochastic systems it is often
necessary to find suitable substitutes.  
In this note we want to draw the attention to an equivalent
formulation of the frequency domain condition, which to our knowledge
is not very present in the literature. We call it the {\em coercivity
  condition}. As our two main contributions, we first establish the equivalence of
the coercivity condition and the frequency domain condition and show
second that the coercivity condition plays the same role for the
solvability of the Riccati equation of a stochastic linear quadratic
control problem as the frequency condition does for the corresponding
deterministic problem. To simplify the presentation we choose the most
simple setup for the stochastic problem. A detailed discussion of the
analogous result for time-varying linear systems is to be found
in the forthcoming book \cite{HinrPrit}. 

It is a great honour for us to dedicate this  note to Vasile
Dr\u{a}gan at the occasion of his 70-th birthday. We had the pleasure to collaborate with Vasile e.g.~in \cite{DragDamm05} and \cite{JacoDrag98}.
Vasile Dr\u{a}gan has made
numerous and substantial contributions in the context of our topic.
Together with Aristide Halanay, he was among the first to study stochastic disturbance attenuation
problems \cite{DragHala96, DragHala97, DragHala99}, and in still
ongoing work (e.g.\ \cite{DragIvan20}) he extended the theory in many different directions. The textbook \cite{DragMoro13} is closely related to this
note.

\section{Preliminaries}
\label{sec:preliminaries}

Consider the time-invariant finite-dimensional linear control system
\begin{eqnarray*}
  \dot x(t)&=&Ax(t)+Bu(t),\quad t\ge 0,\\
  x(0)&=&x_0,
\end{eqnarray*}
together with the {\em  quadratic cost functional}
\begin{displaymath}
  J(x_0,u)=\int_0^\infty \left[
    \begin{array}{c}
      x(t)\\u(t)
    \end{array}
\right]^*M \left[
    \begin{array}{c}
      x(t)\\u(t)
    \end{array}
\right]\,dt\;,
\end{displaymath}
where $A$, $B$ and $M$ are complex matrices of suitable size.
Assume that $M=M^*=\left[
  \begin{array}{cc}
    W&V^*\\V&R 
  \end{array}
\right]$ where $R>0$, but not necessarily $M\ge 0$ or $W\ge 0$. 
With these data we associate the {\em algebraic Riccati equation}
\begin{equation}
  \label{eq:are}
  A^*P+PA+W-(B^*P+V)^*R^{-1}(B^*P+V)=0\;.
\end{equation}
Moreover,  for  $\omega\in{\bf R}$ with 
$\imath\omega\not\in\sigma(A)$  we define the {\em frequency function}
(or Popov function, \cite{Popo73})
\begin{equation}
  \label{eq:freqfunc}
\Phi(\omega)=\left[
    \begin{array}{c}
      (\imath \omega I-A)^{-1}B\\I
    \end{array}
\right]^* M\left[
    \begin{array}{c}
      (\imath \omega I-A)^{-1}B\\I
    \end{array}
\right]\;.
\end{equation}
Here $\sigma(A)$ denotes the spectrum of the matrix $A$.
Then the {\em strict frequency domain condition} requires 
\begin{equation}
  \label{eq:fdc_strict}
  \exists  \varepsilon>0: \forall\omega\in{\bf R},
\imath\omega\not\in\sigma(A):\quad \Phi(\omega)\ge  \varepsilon^2
B ^*(\imath\omega I-A)^* (\imath \omega I-A)^{-1}B\;.
\end{equation}
The {\em nonstrict frequency domain condition} is just 
\begin{equation}
  \label{eq:fdc_nonstrict}
 \forall\omega\in{\bf R},
\imath\omega\not\in\sigma(A):\quad \Phi(\omega)\ge 0\;.
\end{equation}
Note that (\ref{eq:fdc_strict}) holds with a given $M$ and fixed
$ \varepsilon>0$, if and only if (\ref{eq:fdc_nonstrict}) holds with $M$
replaced by
\begin{displaymath}
M_ \varepsilon=M-\left[
  \begin{array}{cc}
     \varepsilon^2I&0\\0&0
  \end{array}
\right]\;.
\end{displaymath}

\begin{remark}\label{rem:fdc_are}
  If $(A,B)$ is stabilizable, then it is well-known (e.g.\
  \cite{Will71}) that (\ref{eq:are}) possesses a {\em stabilizing solution}
  (i.e.\ a solution $P$ with the additional property that
  $\sigma(A-BR^{-1}(B^*P+V))\subset{\bf C}_-$, where ${\bf C}_-$
  denotes the open left half plane), if and only if the frequency
  condition (\ref{eq:fdc_strict}) holds. There exists an {\em almost
    stabilizing solution} (satisfying
  $\sigma(A-BR^{-1}(B^*P+V))\subset{\bf C}_-\cup \imath{\bf R}$), if
  and only if (\ref{eq:fdc_nonstrict}) holds.
\end{remark}

However, there are other classes of linear systems for which quadratic
cost functionals can be formulated, which do not allow for an
analogous frequency domain interpretation.  These are, for instance,
time-varying or stochastic systems e.g.\ \cite{DragMoro13}.
It is therefore useful to have a time domain condition which is
equivalent to (\ref{eq:fdc_strict}). Such a condition can be obtained by
applying the inverse Laplace
transformation, but we choose a more elementary
approach here.

For an initial value $x_0$ and a square-integrable  input function $u\in L^2({\bf R}_+)$ we denote by $x(t,x_0,u)$ the unique solution 
of our time-invariant finite-dimensional linear control system at time $t$.
Let
\begin{displaymath}
U=\{u\in L^2({\bf R}_+)\;\big|\; x(\cdot,0,u)\in L^2({\bf R}_+)\}
\end{displaymath}
denote the set of {\em admissible inputs}. For $u\in U$ we consider 
the cost associated to zero initial state
\begin{equation}
  \label{eq:J0u}
  J(0,u)= \int_0^\infty \left[
    \begin{array}{c}
      x(t,0,u)\\u(t)
    \end{array}
\right]^* M \left[
    \begin{array}{c}
      x(t,0,u)\\u(t)
    \end{array}
\right] \,dt \;.
\end{equation}
Then we say that $J$ satisfies the {\em strict coercivity condition},
if
\begin{equation}
  \label{eq:ccstrict}
  \exists \varepsilon >0 : \forall u\in U: \quad J(0,u)\ge  \varepsilon^2\|x(\cdot,0,u)\|_{L^2}^2\;.
\end{equation}
We say that $J$ satisfies the {\em nonstrict coercivity condition}, if 
\begin{equation}
  \label{eq:cc}
  \forall u\in U:\quad J(0,u)\ge 0\;.
\end{equation}
As for the frequency domain conditions, note that  (\ref{eq:ccstrict}) holds with a given $M$ and fixed
$ \varepsilon>0$, if and only if (\ref{eq:cc}) holds with $M$
replaced by  $M_ \varepsilon$.\\
In the next section, we prove the equivalence of (\ref{eq:fdc_nonstrict}) and
(\ref{eq:cc}). Since in the strict cases with given $ \varepsilon>0$ we
can replace $M$ by $M_ \varepsilon$ as indicated above, this also
establishes the equivalence of  (\ref{eq:fdc_strict}) and
(\ref{eq:ccstrict}). 
Then, in Section \ref{sec:stoch-lq-contr}, we show for a stochastic LQ-problem that
(\ref{eq:ccstrict}) is a natural time domain replacement of (\ref{eq:fdc_strict}).

\section{Equivalence of frequency domain and coercivity condition}
In this section we consider the system $\dot x=Ax+Bu$ and we assume that the pair
$(A,B)\in{\bf C}^{n\times  n}\times {\bf C}^{n\times m}$ is stabilizable.
Solutions with initial value $x(0)=x_0$ and input $u\in
L^2({\bf R}_+)$ are denoted by $x(\cdot,x_0,u)$. As above, let
$$U=\{u\in L^2({\bf R}_+)\;\big|\; x(\cdot,0,u)\in L^2({\bf R}_+)
\}$$  be the set of admissible inputs and let  $M\in{\bf C}^{(n+m)\times(n+m)}$ be a weight matrix of the form $M=M^*=\left[
  \begin{array}{cc}
    W&V^*\\V&R 
  \end{array}
\right]$ where $R>0$.

\begin{theorem} The following statements are equivalent.
  \begin{itemize}
  \item[(a)] For all $u\in U$ it holds that
      \begin{displaymath}
    J(0,u)= \int_0^\infty \left[  \begin{array}{c}
        x(t,0,u)\\u(t)
      \end{array}
    \right]^* M \left[
      \begin{array}{c}
        x(t,0,u)\\u(t)
      \end{array}
    \right] \,dt \ge 0\;.
\end{displaymath}
  \item[(b)] For all $\omega\in{\bf R}$ with
    $\imath\omega\not\in\sigma(A)$  it holds that
    \begin{displaymath}
    \Phi(\omega)=\left[
      \begin{array}{c}
        (\imath \omega I-A)^{-1}B\\I
      \end{array}
    \right]^* M\left[
      \begin{array}{c}
        (\imath \omega I-A)^{-1}B\\I
      \end{array}
    \right] \ge 0\;.
  \end{displaymath}

  \end{itemize}
\end{theorem}
\bprf
(a)$\Rightarrow$(b) Let $\eta\in{\bf C}^m$ be arbitrary and $\omega>0$,
$\imath\omega\not\in\sigma(A)$.  We have to
show that $\eta^*\Phi(\omega)\eta\ge 0$. 
Let $\xi=(\imath\omega I-A)^{-1}B\eta$. Then $\xi$ is reachable from
$0$ and there exists a control input $u_0\in L^2([0,1])$ such that $x(1,0,u_0)=\xi$. Since $(A,B)$ is stabilizable, there also exists $u_\infty\in
L^2({\bf R}_+)$ with $x(\cdot,\xi,u_\infty)\in L^2({\bf R}_+)$.
For $k\in {\bf N}$, $k>0$, and $T_k=\frac{2k\pi}\omega+1 $, we define
\begin{eqnarray*}
  u_k(t)&=&\left\{
          \begin{array}{ll}
            u_0(t)&t\in[0,1[\\
           \eta e^{\imath\omega (t-1)}&t\in[1,T_k]\\
           u_\infty(t-T_k)& t\in\;]T_k,\infty[
          \end{array}
\right.\,.
\end{eqnarray*}
Then $x(1,0,u_k)=\xi$.  An easy calculation shows that on $[1,T_k]$ we have the
resonance solution $x(t,0,u_k)=\xi e^{\imath\omega (t-1)}$ with
$x(T_k,0,u_k)=\xi$, such that it is stabilized by $u_\infty$ on
$]T_k,\infty[$. The integrals
\begin{eqnarray*}
&&\int_0^1 \left[
    \begin{array}{c}
      x(t,0,u_k)\\u_k(t)
    \end{array}
\right]^* M \left[
    \begin{array}{c}
      x(t,0,u_k) \\u_k(t)
    \end{array}
\right] \,dt\\&&+\int_{T_k}^\infty \left[
    \begin{array}{c}
      x(t,0,u_k) \\u_k(t)
    \end{array}
\right]^* M \left[
    \begin{array}{c}
      x(t,0,u_k) \\u_k(t)
    \end{array}
\right] \,dt=c <\infty
\end{eqnarray*}
are independent of  $k$.
By (a) we have
\begin{eqnarray*}
 0&\le& J(0,u_k)=c+\int_1^{T_k} \left[
    \begin{array}{c}
      x(t,0,u_k)\\u_k(t)
    \end{array}
\right]^* M \left[
    \begin{array}{c}
      x(t,0,u_k)\\u_k(t)
    \end{array}
\right] \,dt\\
&=&c+\int_1^{T_k}\left[
    \begin{array}{c}
      \xi e^{\imath\omega t}\\\eta e^{\imath\omega t}
    \end{array}
\right]^* M \left[
    \begin{array}{c}
      \xi e^{\imath\omega t}\\\eta e^{\imath\omega t}
    \end{array}
\right] \,dt\\
&=&c+\int_1^{T_k}\left[
    \begin{array}{c}
      \xi \\\eta
    \end{array}
\right]^* M \left[
    \begin{array}{c}
      \xi \\\eta 
    \end{array}
\right] \,dt\\ &=&c+\int_1^{T_k}\eta^*\left[
    \begin{array}{c}
      (\imath \omega I-A)^{-1}B\\I
    \end{array}
\right]^* M\left[
    \begin{array}{c}
      (\imath \omega I-A)^{-1}B\\I
    \end{array}
\right]\eta \,dt \;.
\end{eqnarray*}
Since $T_k$ can be arbitrarily large, the integrand must be
nonnegative. This proves (b).\\
(b)$\Rightarrow$(a) 
Note first that
\begin{equation}
\left[
  \begin{array}{c}
    \xi\\\eta
  \end{array}
\right]^*M \left[\begin{array}{c}
    \xi\\\eta
  \end{array}
\right]\ge 0, \mbox{ if } 
(\imath\omega I-A) \xi=B\eta \mbox{ for some } \omega\in{\bf R}.\label{eq:initial_observation}
\end{equation}

Let now $u\in U$ be given and assume by way of contradiction that
$J(0,u)<0$.  
For $T>0$, $x_0\in{\bf C}^n$, we set  $$J_T(x_0,u)= \int_0 ^T \left[
    \begin{array}{c}
      x(t,x_0,u)\\u(t)
    \end{array}
\right]^* M \left[
    \begin{array}{c}
      x(t,x_0,u)\\u(t)
    \end{array}
\right] \,dt \;.$$
Then there exists $\delta>0,T_0>1$ such that $J_{T-1}(0,u)<-2\delta$ for
all $T\ge T_0$. For $x_T=x(T-1,0,u)$ there exists a control
input $u_T\in L^2([0,1])$ such that $x(1,x_T,u_T)=0$. 
In fact, one can choose $u_T(t)=-e^{A ^*(1-t)}P_1^\dagger
e^Ax_T$, where $P_1$
denotes the finite-time controllability Gramian over the interval $[0,1]$ and $P_1^\dagger$ its
Moore-Penrose pseudoinverse. Then, e.g.\ \cite{BennDamm11}, 
$$\|u_T\|_{L^2([0,1])}^2=x_T ^*e^{A ^*}P_1^\dagger e^Ax_T={\cal
  O}(\|x_T\|^2)\mbox{ for } x_T\to 0\;. $$ 
This implies that also $J_1(x_{T},u_{T})={\cal O}(\|x_T\|^2)$. Since
$x(\cdot,0,u)\in L^2({\bf R}_+)$, 
we can fix  $T>T_0$ such  that $\|x_{T}\|$  is small enough to ensure
$J_1(x_{T},u_{T})<\delta$. 
We   concatenate $u\big|_{[0,T-1]}$ and $u_T$ to  a new input 
$\tilde u\in L^2([0,T])$ with
\begin{eqnarray*}
  \tilde u(t)&=&\left\{
               \begin{array}{ll}
                u(t) & t\in[0,T-1[\\
               u_T(t-T+1)& t\in[T-1,T]
               \end{array}
\right.\;.
\end{eqnarray*}
By construction, we have 
\begin{equation}
J_T(0,\tilde u)<-\delta<0\quad\mbox{ and }\quad 0=x(0,0,\tilde u) =x(T,0,\tilde
u).\label{eq:JT}
\end{equation}
By definition $\tilde u, x\in L^2([0,T])$, and the equation $\dot x=
Ax +B\tilde u$ implies that $x$ is absolutely continuous and $\dot
x\in  L^2([0,T])$. Thus,  on $[0,T]$, the Fourier series of $\tilde
u$, $x$ and $\dot x$ converge in $ L^2([0,T])$ to $\tilde u$, $x$ and
$\dot x$, respectively.  On $[0,T]$, let 
\begin{displaymath}
x(t,0,u)=\sum_{k=-\infty}^\infty
\xi_ke^{\imath\frac{2\pi k t}T}\quad\mbox{ and }\quad\tilde u(t)=\sum_{k=-\infty}^\infty
\eta_ke^{\imath\frac{2\pi k t}T}.
\end{displaymath}
Then we get
\begin{eqnarray}
\nonumber
\sum_{k=-\infty}^\infty
B\eta_ke^{\imath\frac{2\pi k t}T} =B\tilde u(t) &=&\dot x(t,0,\tilde u)-Ax(t,0,\tilde u)\\&=&\sum_{k=-\infty}^\infty \left(\imath\frac{2\pi}T k
       I-A\right)\xi_ke^{\imath\frac{2\pi k t}T}\;.
\label{eq:formalderivative}
\end{eqnarray}
Note that the periodicity condition $x(0)=x(T)$ in (\ref{eq:JT})
justifies the formal differentiation of the Fourier series in (\ref{eq:formalderivative}),
e.g.\ \cite[Theorem 1]{Tayl44}.\\
Comparing the coefficients in (\ref{eq:formalderivative}), we have
\begin{equation}
  \label{eq:xikBetak}
  \left(\imath\frac{2\pi}T k I-A\right)\xi_k=B\eta_k\;. 
\end{equation}
In the expression of $J_T(0,\tilde u)$ we replace $x(t,0,\tilde u)$ and $\tilde
u(t)$ by
their Fourier-series representations. Exploiting orthogonality we have
\begin{eqnarray*}
J_T(0,\tilde u)&=&
   T\sum_{k=-\infty}^\infty \left[
  \begin{array}{c}
    \xi_k\\\eta_k
  \end{array}
\right]^*M \left[\begin{array}{c}
    \xi_k\\\eta_k
  \end{array}
\right].
\end{eqnarray*}
Together with (\ref{eq:xikBetak}) and (\ref{eq:initial_observation}) this
implies $J_T(0,\tilde u)\ge 0$ contradicting the first condition in
(\ref{eq:JT}). Thus our initial assumption was wrong, and we have shown that
(b) implies (a).\eprf

\section{An indefinite stochastic LQ-control problem}
\label{sec:stoch-lq-contr}

Consider the It\^o-type linear stochastic system
\begin{equation}
  \label{eq:Ito}
  dx=(Ax+Bu)\,dt + Nx\,dw\;.
\end{equation}
Here $w$ is a Wiener process and by $L^2_w$ we denote the space of square
integrable stochastic processes adapted to $L^2_w$. For the
appropriate definitions see textbooks such as \cite{Arno73, DragMoro13}.
Let further the cost functional
\begin{equation}
\label{eq:JWR}
  J(x_0,u)={\bf E}\int_0^\infty \left[
    \begin{array}{c}
      x(t,x_0,u)\\u(t)
    \end{array}
\right]^* M \left[
    \begin{array}{c}
      x(t,x_0,u)\\u(t)
    \end{array}
\right] \,dt 
\end{equation}
be given, where ${\bf E}$ denotes expectation. \\
For simplicity of presentation let $M=\left[
  \begin{array}{cc}
    W&0\\0&I
  \end{array}
\right]$ which can always be achieved by a suitable transformation, if
the lower right block of $M$ is positive definite, e.g.\ \cite[Section
5.1.7]{Damm04}. We do not impose
any definiteness conditions on $W$. Note that we might include further noise processes or control
dependent noise in (\ref{eq:Ito}) at the price of increasing the
technical burden. 

\begin{definition}
  Equation (\ref{eq:Ito}) is {\em internally mean square asymptotically
  stable}, if for all initial conditions $x_0$ the uncontrolled
solution converges to zero in mean square, that is ${\bf
  E}\|x(t,x_0,0)\|^2\to 0$ for $t\to\infty$.
  In this case, for brevity, we also
  call the pair $(A,N)$ {\em asymptotically stable}. We call an input
  signal $u\in L^2_w$ admissible, if also $x(\cdot,0,u)\in L^2_w$.
\end{definition}
It is well known, that the pair $(A,N)$ is asymptotically stable, if
and only if
\begin{displaymath}
    \sigma(I\otimes A+A\otimes I+N\otimes N)\subset{\bf C}_-\;,
  \end{displaymath}
  where $\otimes$ denotes the Kronecker product, \cite{Klei69}.

With (\ref{eq:Ito}) and (\ref{eq:JWR}) we associate the algebraic Riccati equation
\begin{equation}
  \label{eq:Riccati}
  A ^*P+PA+N ^*PN+W-PBB ^*P=0\;.
\end{equation}

\begin{definition}
  A solution $P$ of (\ref{eq:Riccati}) is {\em stabilizing}, if the
  pair $(A-BB ^*P,N)$ is asymptotically stable. We call the triple
  $(A,N,B)$ {\em stabilizable}, if there exists a matrix $F$, such
  that $(A+BF,N)$ is asymptotically stable.
\end{definition}
We now relate the existence of stabilizing solutions of
(\ref{eq:Riccati}) to a coercivity condition. Recall from Remark
\ref{rem:fdc_are} that the frequency condition is used for this
purpose in the deterministic case. For stochastic systems, however,
there is no obvious way to define a transfer function. 
\begin{theorem}
 Let $(A,N,B)$ be {\em stabilizable}.\\
  The Riccati equation (\ref{eq:Riccati}) possesses a stabilizing
  solution, if and only if for some $ \varepsilon>0$ and all admissible
  $u$ the coercivity condition $J(0,u)\ge
   \varepsilon\|x\|_{L_w^2}^2$ holds.  
\end{theorem}
\bprf
We develop the proof along results available in the literature.\\
Let $W=W_1-W_2$, where both $W_1,W_2>0$, and consider first the
definite LQ-problem with the cost functional
\begin{displaymath}
  J_{W_1}(x_0,u)={\bf E}\int_0^\infty (x ^*W_1x+\|u\|^2)\,dt\;.
\end{displaymath}
Then it is known from \cite{Wonh68}, that a minimizing control $u_1$
for $J_{W_1}$ is given in the form $u_1=Fx=-B ^*P_1x$, where $P_1$ is the
unique stabilizing solution of the Riccati equation
\begin{equation}
  \label{eq:Riccati1}
  A ^*P+PA+N ^*PN+W_1-PBB ^*P=0\;.
\end{equation}
For a control of the form $u=-B ^*P_1x+u_2$ it follows that
\begin{equation}\label{eq:J1x0}
   J(x_0,u)=x_0 ^*Px_0+{\bf E}\int_0^\infty \left(\|u_2(t)\|^2-x(t) ^*W_2x(t)\right)\,dt\;,
\end{equation}
where now $x(t)$ is the solution of the closed loop equation
\begin{equation}\label{eq:P1closed}
  dx=(A-BB ^*P_1)x\,dt+Nx\,dw+Bu_2\,dt\;
\end{equation}
with initial value $x_0$. Our next goal is to minimize 
\begin{displaymath}
   J_{W_2}(x_0,u_2)={\bf E}\int_0^\infty \left(\|u_2(t)\|^2-x(t) ^*W_2x(t)\right)\,dt
\end{displaymath}
subject to (\ref{eq:P1closed}).
If we factorize $W_2=C_2 ^*C_2$ and set $y(t)=C_2x(t)$, then we recognize $J_{W_2}$ as the cost
functional related to the stochastic bounded real lemma, \cite[Theorem
2.8]{HinrPrit98}, see also e.g.\ \cite{DragMoro13}.
The associated Riccati {\em inequality}
 \begin{equation}
  \label{eq:Riccati2}
  (A-BB ^*P_1) ^*P+P(A-BB ^*P_1)+N ^*PN+W_2-PBB ^*P>0\;,
\end{equation}
possesses a  solution $\hat P<0$, if and only if there exists a $\delta>0$,
such that
\begin{equation}\label{eq:BRL_condition}
J_{W_2}(0,u_2)\ge\delta\|u_2\|_{L^2_w}^2\mbox{ for all } u_2\in {L^2_w}
\end{equation}
see \cite[Corollary 2.14]{HinrPrit98}. By \cite[Theorem 5.3.1]{Damm04}
this is equivalent to the 
corresponding Riccati {\em equation} having a stabilizing solution
$P_2<0$. \\ 
Note now that for $P=P_1+P_2$, the Riccati equation
(\ref{eq:Riccati}) holds because
\begin{eqnarray*}
  0&=&(A-BB ^*P_1) ^*P_2+P_2(A-BB ^*P_1)+N ^*P_2N+W_2-P_2BB ^*P_2\\ 
&=&A ^*P_2+P_2A+N ^*P_2N-W_2-PBB ^*P+P_1BB ^*P_1\\
&=&A ^*P+PA+N ^*PN+W-PBB ^*P\;.
\end{eqnarray*}
Moreover, the pair $(A-BB ^*P,N)=(A-BB ^*P_1-BB ^*P_2,N)$ is
stabilizing.

It remains to show that (\ref{eq:BRL_condition}) is equivalent to the
coercivity condition. 
As above, let $u\in {L^2_w}$ be of the
form $u=-B ^*P_1x+u_2$. 
Assume first that the coercivity  condition holds. By (\ref{eq:J1x0})  we have
\begin{displaymath}
  J(0,u)=J_{W_2}(0,u_2)=\|u_2\|_{L^2_w}^2-\|y\|_{L^2_w}^2\ge  \varepsilon^2 \|x\|_{L^2_w}^2\;,
\end{displaymath}
where $x$ solves  (\ref{eq:P1closed}) and $y=C_2x$.  It follows that
$\|x\|_{L^2_w}\le \frac1{\|C_2\|}\|y\|_{L^2_w}$, whence
\begin{displaymath}
  \|u_2\|_{L^2_w}^2\ge \left(1+\frac{ \varepsilon^2}{\|C_2\|^2}\right)
  \|y\|_{L^2_w}^2=\alpha  \|y\|_{L^2_w}^2 
\end{displaymath}
with $\alpha>1$. Hence, with $\delta^2=1-\frac1\alpha>0$, we have
\begin{displaymath}
  J_{W_2}(0,u_2)= \|u_2\|_{L^2_w}^2-\|y\|_{L^2_w}^2=\frac1\alpha
  \|u_2\|_{L^2_w}^2-\|y\|_{L^2_w}^2+\delta^2 \|u_2\|_{L^2_w}^2\ge \delta^2 \|u_2\|_{L^2_w}^2 
\end{displaymath}
for all $u_2\in L^2_w$, which is (\ref{eq:BRL_condition}).

Vice versa, assume (\ref{eq:BRL_condition}). 
Since
(\ref{eq:P1closed}) is asymptotically stable, the system has finite
input to state gain $\gamma$, such that $\|x\|_{L^2_w}\le \gamma
\|u_2\|_{L^2_w}$. Therefore 
\begin{displaymath}
  J(0,u)=J_{W_2}(0,u_2)\ge \frac{\delta^2}{\gamma^2} \|x\|_{L^2_w}^2=
  \varepsilon^2\|x\|_{L^2_w}^2
\end{displaymath}
with $ \varepsilon=\frac{\delta}{\gamma}$. This is the coercivity condition.
\eprf

\section{Conclusion}
\label{sec:conclusion}
We have provided a time domain substitute for the frequency domain
condition of the Kalman-Yakubovich-Popov lemma. The equivalence of the
two criteria has been proven and the applicability has been
demonstrated for a stochastic linear quadratic control problem.

\end{document}